\documentclass[11pt]{amsart}
\usepackage{a4,amsfonts,amssymb,amsmath}
\def\proof{{\bf Proof\quad}}

\def\qed{\hfill\rule{2.2mm}{2.2mm}\vspace{1ex}}

\newtheorem{theorem}{Theorem}[section]

\newtheorem{lemma}[theorem]{Lemma}

\def\eps{\varepsilon}

\def\ZZ{\mathbb Z}
\def\DD{\mathcal D}

\def\NN{\mathbb N}
\def\RR{\mathbb R}

\newcommand{\carl}{\mathop{\rm Carl}\nolimits}

\renewcommand{\S}{\mathcal S}
\renewcommand{\span}{\mathop{\rm span}\nolimits}

\newcommand{\diag}{\mathrm{Diag}}
\newcommand{\BMO}{\mathrm{BMO}}
\newcommand{\ohne}{\backslash}

\def\text{\mbox}

\title{Carleson measure and balayage}

\author{Sandra Pott}
\address{Institut f\"ur Mathematik, Universit\"at Paderborn, Warburger Str.~100, 33098 Paderborn \\ and Department of Mathematics, University of Glasgow, Glasgow G12 8QW, UK}
\email{sandrapo@math.uni-paderborn.de}

\author{Alexander Volberg}
\address{Department of Mathematics, Michigan State University, East Lansing, MI 48824, US}
\email{volberg@math.msu.edu}

\begin{document}
\begin{abstract}
The balayage of a Carleson measure lies of course in $\BMO$. We show
that the converse statement is false. We also make a two-sided
estimate of the Carleson norm of a positive measure in terms of
{\em certain} balayages.
\end{abstract}

\maketitle

{\bf Keywords.} BMO functions,  Carleson measure\\

{\bf 2000 Subject Classification.} 42B35

\section{Introduction and notation}

In this note, we consider a question that naturally appeared in the
recent work of Frazier-Nazarov-Verbitsky \cite{fnv}. The question
is:

\medskip

How does the Carleson norm of a positive measure in the disc relate to
the $BMO$ norm of its balayage on the circle?

\medskip

A related question is:

\medskip

How can one describe measures on the disc (say, positive measures)
whose balayage is a $\BMO$-function?

\medskip

The second author is grateful to Igor Verbitsky, who called our
attention to these questions.

We show that the seemingly answer :``These are exactly the Carleson
measures'' is false. The Carleson property is indeed of course
sufficient, but not at all necessary. However, we can characterise the Carleson property
in terms of the $\BMO$ norms of the balayages of restrictions of the measure.

Throughout the paper,
we will use the notation $\lesssim$, $\gtrsim$ for one-sided estimates up to an absolute constant, 
and the notation $\approx$ for two-sided estimates up to an absolute constant.

 We will use the setting of the upper half plane $\RR^2_+$ rather than the unit disc.  Given a positive
regular Borel measure $\mu$ on the upper half plane \\ $\RR^2_+ =
\{(t,y) \in \RR^2: y >0 \}$, its \emph{balayage} is defined as the
function
$$
   S_\mu(t) = \int_{\RR^2_+} p_{x,y}(t) d \mu(t, y),
$$
where $ p_{x,y}(t)= \frac{1}{\pi}\frac{x}{y^2 +(t-x)^2}$ is the Poisson kernel for $\RR^2_+$.
We say that $\mu$ is a \emph{Carleson measure}, if
there exists a constant $C >0$ such that for each interval $I \subset \RR$, the inequality
\begin{equation}  \label{eq:carl}
   \mu(Q_I) \le C |I|
\end{equation}
holds. Here, $Q_I$ denotes the \emph{Carleson square} $\{(x,y): x \in I, 0 < y \le |I|\}$ over $I$. It is easy to see that
it is sufficient to consider dyadic intervals in this definition. We denote the infimum of all constants $C>0$ such that
(\ref{eq:carl}) holds for all dyadic intervals by $\carl(\mu)$.

Recall that the space of functions of \emph{bounded mean oscillation}, $\BMO(\RR)$, is defined as
$$
   \{b \in L^2(\RR): \sup_{I \subset \RR \text{ \small interval}} \frac{1}{|I|} \int_I |b(t) - \langle b \rangle_I|dt < \infty\},
$$
with $\|b\|_{\BMO}=\sup_{I \subset \RR \text{ \small interval}} \frac{1}{|I|} \int_I |b(t) - \langle b \rangle_I|dt$.
By the John-Nirenberg inequality, the $L^1$ norm in the definition of $\BMO$ can be replaced by any $\|\cdot\|_p$
norm, $1\le p < \infty$. We thus obtain a family of equivalent norms on $\BMO(\RR)$, with equivalent
constants depending on $p$.

The connection between the the properties of a measure $\mu$ and
its balayage $S_\mu$ have long been studied. In particular, it is well-known
that the BMO norm of $S_\mu$ is controlled by the Carleson constant of $\mu$,
\begin{equation}  \label{eq:bala}
   \|S_\mu\|_{BMO} \lesssim \carl(\mu).
\end{equation}
For this and other basic facts on
BMO functions, we refer the reader to \cite{garnett}.

A partial reverse of (\ref{eq:bala})
was found in \cite{carl}, \cite{uchi} and, in the dyadic case, \cite{jones}. Namely, it was shown that
for each $b \in \BMO$,
there exists an $L^\infty(\RR)$ function $\phi$ and a Carleson measure $\mu$ such that $b = \phi + S_\mu$,
$\|\phi\|_\infty +\carl(\mu) \lesssim \|b\|_{\BMO}$. If we allow $\mu$ to be a complex measure, one even
has the representation $b = S_\mu$ with $\carl(\mu) \lesssim \|b\|_{\BMO}$ \cite{smith}.

The purpose of this note is to show that reverse inequality to (\ref{eq:bala})
in the strict sense does not hold, and to give a characterization of the
Carleson property of a measure $\mu$ in terms of the BMO norm of the balayage of \emph{restrictions} of $\mu$.

\section{The dyadic balayage}
\label{dyadic}
We start by examining the dyadic case. We will use the standard Whitney-type decomposition
of the upper half-plane, indexed by the set $\DD$ of left-half
open dyadic intervals in $\RR$,
$$
    T_I = \{ (x,y): x \in I, \frac{|I|}{2} < y \le  |I| \} \text{ for } I \in \DD.
$$
That means, $T_I$ is the ``top half'' of the Carleson square $Q_I$ defined above.

For a positive regular Borel measure $\mu$ on $\RR^2_+$, we define the \emph{dyadic balayage} by
$$
   S^d_\mu(t) = \sum_{I \in \DD} \frac{\chi_I(t)}{|I|} \mu(T_I)  \quad (t \in \RR),
$$
which is well-defined as a function taking values in $[0,\infty]$. By comparing box kernel and Poisson kernel,
one easily verifies the pointwise estimate $S^d_\mu \lesssim S_\mu$.

We recall the definition of \emph{dyadic BMO}, $\BMO^d(\RR)$, as the class of $L^2(\RR)$ functions for which
$$
   \|b\|_{BMO^d}^2 = \sup_{I \in \DD} \frac{1}{|I|} \int_I |b(t) - \langle b \rangle_I|^2 dt =
 \sup_{I \in \DD}\frac{1}{|I|} \|P_I b\|^2 =  \sup_{I \in \DD}\frac{1}{|I|} \sum_{J \in \DD, J \subseteq I} |b_I|^2
$$
is finite. Here, $h_J$ denotes the $L^2$-normalized Haar function,
$b_J:=(b, h_J)$ denotes the corresponding Haar coefficient of
function $b$,  and $P_I$ denotes the orthogonal projection onto
$\overline{\span \{h_J: J \subseteq I\}}$. Again, by the
John-Nirenberg inequality the $L^2$ norm in the definition can be
replaced by any $L^p$ norm, $1\le p < \infty$, yielding an
equivalent norm.

We say that a sequence of nonnegative numbers $(\alpha_I)_{I \in \DD}$ is a \emph{Carleson sequence}, if
there exists a constant $C>0$ such that
$$
    \frac{1}{|I|} \sum_{J \in \DD, J \subseteq I} a_I \le C \text{ for each } I \in \DD.
$$
Again, we denote the infimum of such constants by $\carl((a_I))$.
 With this notation, one verifies immediately the following well-known lemma.
\begin{lemma} \label{lemm:carlseq} Let $b \in L^2(\RR)$. Then the following are equivalent:
\begin{enumerate}
\item $\mu$ is a Carleson measure
\item $(\mu(T_I))_{I \in \DD}$ is a Carleson sequence
\item $b_\mu = \sum_{I \in \DD} h_I \mu(T_I)^{1/2} \in \BMO^d(\RR)$.
\end{enumerate}
In this case, $\carl(\mu) = \carl((\mu(T_I))) = \|b_\mu\|^2_{\BMO^d}$.
\end{lemma}

Notice that with the above definition of $b_\mu$,
$$
   S^d_\mu = \sum_{I \in \DD} \frac{\chi_I}{|I|} \mu(T_I)  = \sum_{I \in \DD} \frac{\chi_I}{|I|} |(b\mu)_I|^2
=\S[{b_\mu}],
$$
where $\S$ denotes the square of the dyadic square function, $\S[f] = \sum_{I \in \DD} \frac{\chi_I}{|I|} |f_I|^2$
for $f \in L^2(\RR)$. In this sense, we have identified the dyadic balayage of a positive regular Borel measure
$\mu$ with the square of a dyadic square function of $b_\mu$.
Conversely, for any $f \in L^2(\RR)$, $\S[f]$ can be written as a dyadic balayage of a measure $\mu_f$,
for example by letting $\mu_f = \sum_{I \in \DD} |f_I|^2 \delta_{z(I)}$, $z(I)$ denoting the center of $T_I$.

The well-known dyadic analogue of (\ref{eq:bala}) is therefore equivalent
to the inequality
\begin{equation} \label{eq:dbala}
   \|\S[b]\|_{\BMO^d} \lesssim \|b\|_{\BMO^d}^2,
\end{equation}
which can be now be
proved as a simple application of the John-Nirenberg inequality. Notice that for any dyadic inverval
$I \in \DD$, all summands in $\S[b] = \sum_{J \in \DD} \frac{\chi_J}{|J|} |b_J|^2$ except those
corresponding to dyadic intervals $J \subset I$ are  constant on $I$. Thus
\begin{multline*}
   \frac{1}{|I|}\int_I | \S[b](t) - \langle \S[b] \rangle_I | dt
   = \frac{1}{|I|}\int_I | \S[P_I b](t) - \langle \S[P_Ib] \rangle_I | dt \\
  \le  \frac{1}{|I|}\int_I  \S[P_I b](t) dt + \langle \S[P_Ib] \rangle_I
= 2 \frac{1}{|I|} \int_I \sum_{J \subseteq I} \frac{\chi_J(t)}{|J|} |b_J|^2 dt
=2 \|P_I b\|_2^2 \le 2 \|b\|^2_{\BMO^d},
\end{multline*}
which proves (\ref{eq:dbala}).
\qed

Here are the main results of this section, which concern the reverse inequality to
(\ref{eq:dbala}). The first says that the $\BMO$ norm of the dyadic balayage can be
very much smaller than the Carleson constant of a measure, even if one increases the
$\BMO$ norm by the $L^2$ norm.
\begin{theorem} \label{thm:dcounter} Let $\eps >0$. Then there exists a Carleson measure $\mu$ on $\RR^2_+$ with
$\carl(\mu)=1$, $\|S^d_\mu\|_{\BMO} + \|S^d_\mu\|_2 < \eps$.
\end{theorem}
\proof By Lemma \ref{lemm:carlseq} and the argument following it,
we want to find a $\BMO^d(\RR)$ function $b$ of norm $1$ such that both the $\BMO^d$ norm and the $L^2$ norm
of $\S[b]$ are small. To this end, let $I_{0}= (0, 1]$, $I_{-1}=(-2,0]$, $I_k = (2^{k}-1, 2^{k+1}-1]$ for $k > 0$ and
$I_{k}= (-2^{-k},-2^{-k-1}]$ for $k<0$. In particular, $|I_k| = 2^{|k|}$ for all $k \in \NN$.
Let $r_1$ denote the first Rademacher function on $\RR$, $r_1 =
\sum_{j \in \ZZ} (-\chi_{(j, j+ \frac{1}{2}]} +\chi_{(j+ \frac{1}{2}, j+1]})$, and let
$r_n = r_1(2^{n-1} \cdot)$ be the $n$th Rademacher function on $\RR$. Let $N \in \NN$, $N$ to be determined later, and let
$$
  b= \sum_{k=-\infty}^\infty \sum_{n=1}^{N-|k|} \chi_{I_k}(t) r_n(t).
$$
One verifies without difficulty that $\|b\|_{\BMO^d}^2 =N$. Clearly
$$
  \S[b] = \sum_{k=-\infty}^\infty \sum_{n=1}^{N-|k|} \chi_{I_k} = \sum_{k=0}^N (N-k) \chi_{I_k \cup I_{-k}}.
$$
This is a ``dyadic $\log$'', and it is not difficult to show that
$$
   \|\S[b]\|_{\BMO} \le C,
$$
where $C$ is an absolute constant independent of $N$.  Notice that we have an estimate here
not only for the dyadic $\BMO$ norm, but for the full $\BMO$ norm.

Now choose $N$ so large that $\frac{C}{N} < \frac{\eps}{2}$ and
  replace $b$ by $\frac{1}{N^{1/2}} b$. This already guarantees that
$\|b\|_{\BMO^d}^2 =1$, $\|\S[b]\|_{\BMO}< \frac{\eps}{2}$.
To deal with the desired $L^2$ estimate, observe that
the estimates achieved so far do not change at all if $b$ is dilated with an integer power of $2$.
By choosing a suitable power $2^K$ of $2$, $K \in \NN$, and replacing $b$ by $b(2^K \cdot)$, we obtain the
desired estimate
$$
\|b\|_{\BMO^d}^2 =1,  \quad \|\S[b]\|_{\BMO} + \|S[b]\|_2 < \eps.
$$
\qed

The next theorem says that we can retrieve the Carleson constant of a measure up to an absolute constant  from its
dyadic balayage, if we restrict the measure to certain sets. 
\begin{theorem} \label{thm:dbalay} Let $\mu$ be Carleson measure $\mu$ on $\RR^2_+$. Then
$$
 \carl(\mu) \approx \sup_{E \subseteq \RR^2_+, E \text{ Borel set}} \|S^d_{\mu_E}\|_{\BMO^d}
\approx  \sup_{I \in \DD} \|S^d_{\mu_{Q_I}}\|_{\BMO^d}.
$$
Here, $\mu_E$ stands for the restriction of $\mu$ to $E$, given by $\mu_E(A) = \mu(E \cap A)$.
\end{theorem}
\proof
Clearly $\carl(\mu_E) \le \carl(\mu)$ for each Borel set $E \subseteq \RR^2_+$, so
\begin{multline*}
   \sup_{I \in \DD} \|S^d_{\mu_{Q_I}}\|_{\BMO^d} \le
     \sup_{E \subseteq \RR^2_+, E \text{ Borel set}} \|S^d_{\mu_E}\|_{\BMO^d} \\
   \lesssim \sup_{E \subseteq \RR^2_+, E \text{ Borel set}} \carl(\mu_E) \le \carl(\mu).
\end{multline*}
To prove the reverse inequality, let $I \in \DD$. Observe that $S^d_{\mu_{Q_I}}$ is supported on the closure
of $I$. Therefore, with $I'$ denoting the dyadic sibling of $I$, we have
\begin{multline*}
 \|S^d_{\mu_{Q_I}}\|_{\BMO^d} \ge |\langle S^d_{\mu_{Q_I}}\rangle_I -\langle S^d_{\mu_{Q_I}}\rangle_{I'}|
=\langle S^d_{\mu_{Q_I}}\rangle_I  \\ =
   \frac{1}{|I|} \int_I \sum_{J \in \DD, J \subseteq I} \frac{\chi_J(t)}{|J|} \mu(T_J) dt = \frac{1}{|I|}\mu(Q_I).
\end{multline*}
Thus
$\carl(\mu) \lesssim  \sup_{I \in \DD} \|S^d_{\mu_{Q_I}}\|_{\BMO^d}$.
\qed

\section{The Algebra of Paraproducts}
\label{para}
This section contains a short operator-theoretic motivation for the choice of the counterexample, in
particular the  appearance of Rademacher functions, in the previous section, in terms of \emph{paraproducts}.
Recall that for $b \in L^2(\RR)$, the standard dyadic paraproduct $\pi_b$ is defined by
$$
   \pi_b f = \sum_{I \in \DD} h_I  b_I  \langle f \rangle_I  \text{ for } f \in L^\infty(\RR) \cap L^2(\RR).
$$ 
It is well known, and indeed a reformulation of the classical Carleson Embedding Theorem, that $\pi_b$ extends to a bounded linear operator on $L^2(\RR)$,
if and only if $b \in \BMO^d(\RR)$. In this case, $\| \pi_b\| \approx \|b\|_{\BMO^d}$.

Such dyadic paraproducts have the nice property that $\pi_b^* \pi_b$ is essentially a dyadic paraproduct again, with symbol $\S[b] $ (see \cite{bp}):
\begin{equation} \label{eq:paraid}
   \pi_b^* \pi_b = \pi_{\S[b]} + (\pi_{\S[b]})^* + \diag(b),
\end{equation}
where $\diag(b)$ denotes the diagonal of $\pi_b^*\pi_b$ with respect to the the Haar basis, $\diag(b) h_I = \| \pi_b h_I \|^2 h_I$ for $I \in \DD$.
Moreover,
\begin{equation}   \label{eq:norms}
\| \pi_{\S[b]} \| \approx  \| \pi_{\S[b]} + (\pi_{\S[b]})^* \| \approx \|S[b] \|_{\BMO^d} .
\end{equation}
As pointed out in the previous section, the problem of finding a Carleson measure with Carleson
constant $1$ and small $\BMO^d$ norm of the dyadic balayage is equivalent to finding  $b \in \BMO^d(\RR)$ of norm $1$ such that $\S[b]$
has small $\BMO^d$ norm.

In  light of (\ref{eq:paraid}) and (\ref{eq:norms}), this means finding $b \in \BMO^d(\RR)$ such that $\pi_b^* \pi_b$ is   ``almost diagonal", in the
sense that 
$$
    \|\S[b]\|_{\BMO^d}  \approx \|\pi_{\S[b]} + (\pi_{\S[b]})^* \|=  \|\pi_b^* \pi_b - \diag_b \| \ll  \|\pi_b^*\pi_b\|= \|\pi_b\|^2 \approx \|b\|_{\BMO^d}^2.
$$

Note the elementary identity 
\begin{equation}   \label{eq:diagpart}
   \pi_b^* \pi_b h_I = \frac{1}{|I|^{1/2}} \left( \sum_{J \subseteq I^+}  \frac{\chi_J}{|J|} |b_J|^2 -\sum_{J \subseteq I^-}  \frac{\chi_J}{|J|} |b_J|^2 \right).
\end{equation}
The function $\sum_{J \subseteq I^+}  \frac{\chi_J}{|J|} |b_J|^2  + \sum_{J \subseteq I^-}  \frac{\chi_J}{|J|} |b_J|^2$ is constant on its support $I$ for each $I$,
 if $b$ is a sum of Rademacher functions. In this case, the right-hand side
$\sum_{J \subseteq I^+}  \frac{\chi_J}{|J|} |b_J|^2 -\sum_{J \subseteq I^-}  \frac{\chi_J}{|J|} |b_J|^2$  of (\ref{eq:diagpart}) is always a multiple of $h_I$, and $\pi_b^* \pi_b$ is diagonal in the Haar basis.
In our counterexample, we have to introduce cutoffs on the Rademacher functions in order to control the $L^2$ norm. This introduces nondiagonal terms,
but these can then be controlled by the logarithmic staggering of the cutoffs.

\section{The Poisson balayage}
\label{poisson}

We are now going to construct a compactly supported positive measure $\mu$ on the upper half-plane such that its Carleson constant $\carl(\mu)$ is very
large (say $m$), but $\|S_{\mu}\|_{\BMO} + \|S_{\mu}\|_{L^1}$ is
bounded by absolute constant. From here one can
easily construct finite positive measure $\mu$ which is not
Carleson, but whose balayage is a nice $\BMO$ function.

Fix $m \in \NN$. For $0 \le j \le m$, let $I_j$ denote the interval $[-2^{j}, 2^{j}]$ and
$\tilde I_j = I_j \ohne I_{j-1}$. Furthermore, let $ \tilde I_0 = I_0$ and let $\tilde I_{m+1} = \RR \ohne I_m$.

Let $\mu_j$ denote one-dimensional Lebesgue measure
on the segment $I_j \times \{2^{-j}\}$, and let $\mu = \sum_{j=0}^m m_j$. Clearly $\carl(\mu) =m+1$.

Here is the elementary technical lemma which will show the desired properties of $\mu$. 
\begin{lemma}
There exists an absolute constant $c>0$ (independent of $m$)  such that 
$$
   |S_{\mu_j}(t) - \chi_{I_j}(t)| \le  c \, 2^{-2j}   \text{ for } |t|  \le 2^{j-1}  \text{ or  } |t| \ge 2^{j+1}, \; j \in \{0, \dots, m\}.
$$

\end{lemma}
\proof
Observe that
\begin{eqnarray*}
S_{\mu_j}(t) &= &\frac{1}{\pi} \int_{-2^{j}}^{2^j}  \frac{2^{-j}}{ (x-t)^2 + 2^{-2j}} dx  
                          \le S_{\mu_j}(0)  \le 1 \text{ for all } t \in \RR, \; j \in \{0, \dots, m\}. 
\end{eqnarray*}
Now let $|t| \le 2^{j-1}$. Then
\begin{eqnarray*}
     S_{\mu_j}(t) - 1 &=&\frac{1}{\pi} \int_{-2^{j}}^{2^j}  \frac{2^{-j}}{ (x-t)^2 + 2^{-2j}} dx  - \frac{1}{\pi} \int_{-\infty}^{\infty}  \frac{2^{-j}}{ (x-t)^2 + 2^{-2j}} dx \\
     &=&\frac{1}{\pi} \int_{-\infty}^{-2^j}  \frac{2^{-j}}{ (x-t)^2 + 2^{-2j}} dx+ \frac{1}{\pi} \int_{2^j}^{\infty}  \frac{2^{-j}}{ (x-t)^2 + 2^{-2j}} dx\\
     &\le& \frac{2}{\pi} \int_{0}^{\infty}  \frac{2^{-j}}{ (x+2^{j-1})^2 + 2^{-2j}} dx \\
     &=&  \frac{2}{\pi} \int_{2^{2j-1}}^{\infty}  \frac{1}{ x^2 + 1} dx  \le  
           \sum_{l=j}^\infty   \frac{2}{\pi}  \int_{2^{2l-1}}^{2^{2l+1}} \frac{1}{x^2+1} dx \\
                                                      &\le& \frac{6}{\pi} \sum_{l=j}^\infty 2^{2l-1}   \frac{1}{(2^{2l-1} )^2} = \frac{8}{\pi} 2^{-2j+1}.
\end{eqnarray*}

If $|t| \ge 2^{j+1}$, then
\begin{eqnarray*}
   S_{\mu_j}(t) &= &\frac{1}{\pi} \int_{-2^{j}}^{2^j}  \frac{2^{-j}}{ (x-t)^2 + 2^{-2j}} dx  \\
                         &\le &\frac{1}{\pi} \int_{-2^{j}}^{2^j}  \frac{2^{-j}}{ 2^{2j} + 2^{-2j}} dx  \\
                         & \le & \frac{1}{\pi} 2^{-2j+1}.
\end{eqnarray*}
\qed

Writing $S_\mu = \sum_{j=0}^m S_{\mu_j} = \sum_{j=0}^m \chi_{I_j} + \sum_{j=0}^m ( S_{\mu_j}- \chi_{I_j})$, we see that the first term is a dyadic log function, and therefore in $\BMO(\RR)$ with some
absolute norm bound independent of $m$.  To estimate the second term,  let $t \in \tilde I_k$. By the previous lemma, 
$|S_{\mu_j}(t)- \chi_{I_j}(t)| \le c \,2^{-j}$ for $j \notin \{k-1, k, k+1\}$, therefore
$$
    \sum_{j=0}^m |S_{\mu_j}(t)- \chi_{I_j}(t)|  \le  \sum_{j=0}^m c\, 2^{-j}   + 6 = 2c+6.
$$
Thus the second term is in $L^\infty(\RR)$, with $L^\infty$ norm bounded by $2c +6$. Altogether, we find that there is an absolute
constant $\tilde c$, independent of $m$, such that
$
        \|S_\mu\|_{\BMO} \le \tilde c
$.
However, an elementary calculation shows that 
$$
\|S_{\mu}\|_1 = \sum_{j=0}^m \|S_{\mu_j}\|_1 = \sum_{j=0}^m 2^{j+1}= 2^{m+2}-2,
$$
 and we would like to control the $L^1$ norm of $S_\mu$
as well.
But by scaling our construction with a small $h>0$, i.e.~replacing each $\mu_j$ by $\tilde \mu_j$, the one-dimensional Lebesgue measure on $[- h 2^{j}, h 2^{j}] \times \{h 2^{-j}\}$
and letting $\tilde \mu = \sum_{j=0}^m \tilde \mu_j$, we obtain a measure $\tilde \mu$ with
$\carl(\tilde \mu)=\carl(\mu)= m+1 $,
 $S_{\tilde \mu}(t) = S_\mu(\frac{t}{h})$. Thus  we have $\|S_{\mu}\|_1 = h (2^{m+2}-2)$ and  $\|S_{\tilde \mu} \|_{\BMO} = \|S_{\mu}\|_{\BMO} \le \tilde c$.
 
 After choosing an appropriate $h >0$ and dividing by an appropriate multiple of $m$, we obtain
\begin{theorem} \label{thm:pcounter} Let $\eps >0$. Then there exists a Carleson measure $\mu$ on $\RR^2_+$ with
$\carl(\mu)=1$, $\|S_\mu\|_{\BMO} + \|S_\mu\|_1 < \eps$.
\end{theorem}

We will now show a continuous analogue to Theorem
 \ref{thm:dbalay}.
 \begin{theorem} \label{thm:balay}
 Let $\mu$ be Carleson measure $\mu$ on $\RR^2_+$. Then
$$
 \carl(\mu) \approx \sup_{E \subseteq \RR^2_+, E \text{  \small Borel set}} \|S^d_{\mu_E}\|_{\BMO^d}
\approx  \sup_{I \subset \RR \text{ \small interval}} \|S_{\mu_{Q_I}}\|_{\BMO}.
$$
\end{theorem}

\proof
We only have to prove that 
$\sup_{I \subset \RR \text{ \small interval}} \|S{\mu_{Q_I}}\|_{\BMO} \gtrsim  \carl(\mu)$. After translation and dilation of $\mu$, we can assume without loss of generality that
$\mu(Q_{J}) \ge \frac{1}{4} \carl(\mu)$ for $J= [1/4,3/4]$. Let $I=[0,1]$ and let  $I' $ denote the translated interval $[2,3]$. 
Then
\begin{eqnarray*}
   \|S_{\mu_{Q_I}}\|_{\BMO} &\gtrsim&  |\langle S_{\mu_{Q_I}} \rangle_I - \langle S_{\mu_{Q_I}} \rangle_{I'} | \\
   &=& \int_0^1\frac{1}{\pi} \int_{Q_I}  \frac{y}{(t-x)^2 +y^2}  -  \frac{y}{(t+2 -x)^2 +y^2}d\mu(x,y) dt  \\
   &=& \frac{1}{\pi} \int_{Q_I} \int_{-x}^{1-x} \frac{y(4+4t)}{ (t^2 +y^2)( (t+2)^2 +y^2)} dt d\mu(x,y)\\
    &\ge & \frac{1}{\pi} \int_{[1/4, 3/4] \times [0,1]}   \int_{-x}^{1-x} \frac{y(4+4t)}{ (t^2 +y^2)( (t+2)^2 +y^2)} dt d\mu(x,y)\\
      &\ge & \frac{1}{\pi} \int_{[1/4, 3/4] \times [0,1]}   \int_{-1/4}^{1/4} \frac{y(4+4t)}{ (t^2 +y^2)( (t+2)^2 +y^2)} dt d\mu(x,y)\\
        &\gtrsim & \frac{1}{\pi} \int_{[1/4, 3/4] \times [0,1]}   \int_{-1/4}^{1/4} \frac{y}{ t^2 +y^2} dt\,  d\mu(x,y)\\
        &\ge& \frac{1}{\pi} \int_{[1/4, 3/4] \times [0,1]}   \int_{-1/4}^{1/4} \frac{1}{ t^2 +1} dt \,  d\mu(x,y) \gtrsim \mu(Q_J) \gtrsim  \carl(\mu).\\
\end{eqnarray*}

\qed

\end{document}